\documentclass[11pt]{article}

\usepackage{amssymb,latexsym}
\usepackage{amsmath,amscd}
\usepackage{citesort}
\usepackage{theorem}
\usepackage{eucal}
\usepackage{bbm}
\usepackage[latin1]{inputenc}

\title{The Picard groupoid in deformation quantization\thanks{Talk given at the PQR 2003 Euro-Conference in June 2003, Brussels.}\addtocounter{footnote}{5}}

\author{\textbf{Stefan Waldmann}\thanks{E-mail: Stefan.Waldmann@physik.uni-freiburg.de}
  \\[0.1cm]
  Fakult{\"a}t f{\"u}r Mathematik und Physik\\
  Albert-Ludwigs-Universit{\"a}t Freiburg\\
  Physikalisches Institut\\
  Hermann Herder Stra{\ss}e 3\\
  D 79104 Freiburg\\
  Germany}

\date{November 2003\\[0.5cm] FR-THEP 2003/22
\\[1cm]
\emph{Dedicated to Alan Weinstein for his 60th birthday}}

\renewcommand{\mathbb}[1]{\mathbbm{#1}} 

\newcommand{\im}         {\mathrm{i}}
\newcommand{\cc} [1]     {\overline{{#1}}}
\newcommand{\id}         {{\mathsf{id}}}
\newcommand{\Hom}        {{\mathrm{Hom}}}
\newcommand{\SP} [1]     {{\left\langle{{#1}}\right\rangle}}
\newcommand{\cl}         {{\mathrm{cl}}}
\newcommand{\Def}        {{\mathrm{Def}}}
\newcommand{\ring}[1]    {{\mathsf{{#1}}}}
\newcommand{\Bimod}[3] {{\sideset{_{\scriptscriptstyle{#1}}}{_{\scriptscriptstyle{#3}}}{\operatorname{#2}}}}
\newcommand{\BEA}    {\Bimod{\mathcal{B}}{\mathcal{E}}{\mathcal{A}}}   
\newcommand{\AAA}   {\Bimod{\mathcal{A}}{\mathcal{A}}{\mathcal{A}}}   

\newcommand{\DED}   {\Bimod{\mathcal{D}}{\mathcal{E}}{\mathcal{D}}}
\newcommand{\CFB}   {\Bimod{\mathcal{C}}{\mathcal{F}}{\mathcal{B}}}

\newcommand{\HD}    {\Bimod{}{\mathcal{H}}{\mathcal{D}}}
\newcommand{\HpA}   {\Bimod{}{\mathcal{H}^\prime}{\mathcal{A}}}
\newcommand{\HB}    {\Bimod{}{\mathcal{H}}{\mathcal{B}}}
\newcommand{\HDp}   {\Bimod{}{\mathcal{H}^\prime}{\mathcal{D}}}
\newcommand{\DccH}  {\Bimod{\mathcal{D}}{\cc{\mathcal{H}}}{}}

\newcommand{\tensor}[1][{}]{\mathbin{\otimes_{\scriptscriptstyle{\mathcal{{#1}}}}}}
\newcommand{\tensM}[1][{}] {\mathbin{\widehat{\otimes}_{\scriptscriptstyle{\mathcal{{#1}}}}}}
\newcommand{\rep}[1][{}]  {\sideset{^*}{_{\mathcal{#1}}}{\operatorname{\textrm{-}\mathsf{rep}}}}
\newcommand{\Rep}[1][{}]  {\sideset{^*}{_{\mathcal{#1}}}{\operatorname{\textrm{-}\mathsf{Rep}}}}
\newcommand{\smod}[1][{}]  {\sideset{^*}{_{\mathcal{#1}}}{\operatorname{\textrm{-}\mathsf{mod}}}}
\newcommand{\Pic}  {\mathsf{Pic}}
\newcommand{\StrPic} {\mathsf{Pic}^{\mathrm{str}}}
\newcommand{\starPic} {\mathsf{Pic}^*}
\newcommand{\IP}[4]{{\,}_{\scriptscriptstyle{#2}\!\!}\left\langle{{#1}}\right\rangle^{\scriptscriptstyle{#3}}_{\scriptscriptstyle{#4}}}
\newcommand{\SPEA}[1] {\IP{{#1}}{}{\mathcal{E}}{\mathcal{A}}}
\newcommand{\SPFB}[1] {\IP{{#1}}{}{\mathcal{F}}{\mathcal{B}}}
\newcommand{\SPFEA}[1] {\IP{{#1}}{}{\mathcal{F}\otimes\mathcal{E}}{\mathcal{A}}}
\newcommand{\SPA}[1] {\IP{{#1}}{}{}{\mathcal{A}}}
\newcommand{\ASP}[1] {\IP{{#1}}{\mathcal{A}}{}{}}
\newcommand{\BSP}[1] {\IP{{#1}}{\mathcal{B}}{}{}}
\newcommand{\SPD}[1] {\IP{{#1}}{}{}{\mathcal{D}}}

\newtheorem{proposition}{Proposition}[section]
\newenvironment{proof}[1][{}]{\par\noindent\textsc{Proof{#1}: }}{\hspace*{\fill}$\blacksquare$\smallskip\noindent\par}

\numberwithin{equation}{section}

\begin{document}

\maketitle

\begin{abstract}
    In this letter we give an overview on recent developments in
    representation theory of star product algebras. In particular, we
    relate the $^*$-representation theory of $^*$-algebras over
    rings $\ring{C} = \ring{R}(\im)$ with an ordered ring $\ring{R}$
    and $\im^2 = -1$ to the $^*$-representation theory of
    $^*$-algebras over $\mathbb{C}$ and point out some properties of
    the Picard groupoid corresponding to the notion of strong Morita
    equivalence. Some Morita invariants are interpreted as arising
    from actions of this groupoid. 
\end{abstract}

%
%

\section{Introduction}
\label{sec:intro}

The purpose of this letter is to review some recent developments in
deformation quantization \cite{bayen.et.al:1978a} linked to the
question of finding and describing a physically useful representation
theory for the star product algebras. In fact, this question will be
embedded in a larger context of representation theory of associative
algebras, which arises in many flavours illustrated by the following
(certainly incomplete) list:
\begin{center}
    \begin{tabular}{ccccccccc}
        \multicolumn{9}{c}{Representation theory of} \\[0.2cm]
        rings 
        & $\supseteq$ 
        & rings with 
        & $\supseteq$ 
        & $^*$-algebras
        & $\supseteq$ 
        & $^*$-algebras 
        & $\supseteq$ 
        & $C^*$/$W^*$-algebras \\
        &
        & involutions
        &
        & over $\ring{C} = \ring{R}(\im)$
        &
        & over $\mathbb{C}$
    \end{tabular}
\end{center}
Here and in the following $\ring{R}$ denotes an ordered ring and
$\im^2 = -1$. Star products belong to the middle class of algebras,
whence we shall mainly focus on this type. However, it is clear that
the remaining classes are of great importance in various areas of
mathematics and physics as well.

One main theme of our point of view is that to each version of
representation theory there is (or at least should be) an appropriate
notion of Morita equivalence of the underlying algebras implying in
particular that Morita equivalent algebras have equivalent categories
of representations, see e.g.
\cite{morita:1958a,rieffel:1974b,ara:1999a,bursztyn.waldmann:2003a:pre}
and references therein for the corresponding notions of Morita
equivalence.  Furthermore, Morita equivalence is encoded in the
existence of particular bimodules and collecting all such `equivalence
bimodules' modulo bimodule isomorphisms gives a (large) groupoid, the
Picard groupoid in each of the above situations, where the
multiplication is induced by tensor products of the bimodules.
Moreover, the `equivalence bimodules' implement the functorial
equivalences of representation theories by tensoring with them. Thus
one finally arrives at the picture that the Picard groupoid `acts' on
the representation theories. In this letter we try to make these ideas
more precise, in particular for the case of $^*$-algebras over
$\ring{C} = \ring{R}(\im)$.

Beside this algebraic framework of Morita theory there are many other
notions of Morita equivalence in different areas of
mathematics. Important for us is Xu's notion of Morita equivalent
Poisson manifolds \cite{xu:1991a} which in some sense should be the
`classical limit' of the Morita equivalence of star product
algebras. However, the precise relation is still to be explored, see
in particular the discussion in
\cite{bursztyn.weinstein:2003a:pre,bursztyn.radko:2003a}.

The paper is organized as follows: In Section~\ref{sec:StarToStar} we
review the motivation why one should take a look at
$^*$-representations of star product algebras. In
Section~\ref{sec:positiv} we discuss some general features of
$^*$-representation theory, in particular the notions of positivity
arising from the ordered ring $\ring{R}$. Here we clarify some
relations to the $^*$-representation theory of complex $^*$-algebras.
Section~\ref{sec:TensorProd} is devoted to pre-Hilbert modules and
their tensor products and gives another interpretation of the complete
positivity of inner products. In Section~\ref{sec:Morita} we discuss
strong Morita equivalence and the resulting strong Picard groupoid
while in the last section we demonstrate how some Morita invariants
can be seen as arising from actions of the Picard groupoid.

\medskip

\noindent
\textbf{Acknowledgements:} It is a pleasure to thank the organizers and
in particular Simone Gutt for the excellent summer school and
conference. Furthermore, I would like to thank Mohamed Barakat,
Henrique Bursztyn, Prof.~Schmüdgen, and Alan Weinstein for valuable
discussion and remarks.

%
%

\section{From star products to $^*$-representation theory}
\label{sec:StarToStar}

Let $(M, \pi)$ be a Poisson manifold. Then a formal star product
$\star$ is a $\mathbb{C}[[\lambda]]$-bilinear associative product for
$C^\infty(M)[[\lambda]]$, written as
\begin{equation}
    \label{eq:Starproduct}
    f \star g = \sum_{r=0}^\infty \lambda^r C_r(f, g),
\end{equation}
such that $C_0(f,g) = fg$ is the pointwise `classical' product and
$C_1(f, g) - C_1(g,f) = \im \{f, g\}$ gives the Poisson bracket
induced by $\pi$. Furthermore, one requires $1 \star f = f = f \star
1$ and the $C_r$ are bidifferential operators
\cite{bayen.et.al:1978a}. If $\star$ satisfies in addition
\begin{equation}
    \label{eq:ccfstarg}
    \cc{f \star g} = \cc{g} \star \cc{f},
\end{equation}
where $\lambda = \cc{\lambda}$, then $\star$ is called a Hermitian
star product. Two star products $\star$ and $\star'$ are called
equivalent if there exists an operator $S = \id + \sum_{r=1}^\infty
\lambda^r S_r$ where $S_r$ is a differential operator with $S_r1 = 0$
such that $S(f \star g) = Sf \star' Sg$.  Physically interpreted, the
algebra $(C^\infty(M)[[\lambda]], \star)$ is the algebra of
observables of the quantum theory corresponding to the classical
mechanical system described by $(M, \pi)$ and $\lambda$ plays the role
of Planck's constant $\hbar$. The existence of star products and their
classification up to equivalence is by now well-understood, both for
the symplectic case
\cite{dewilde.lecomte:1983b,fedosov:1994a,omori.maeda.yoshioka:1991a,bertelson.cahen.gutt:1997a,nest.tsygan:1995a,weinstein.xu:1998a}
and for the more general Poisson case
\cite{kontsevich:1997:pre,cattaneo.felder:2000a}. See e.g.
\cite{dito.sternheimer:2002a,gutt:2000a} for recent reviews.

In order to have a physically meaningful quantization one needs more
than the algebra of observables: also the states have to be described.
Since in deformation quantization the observable algebra is realized
as the primary object one seeks for a derived description for the
states. Similarly to and in fact motivated by algebraic quantum field
theory, see e.g.~\cite{haag:1993a}, the \emph{positive linear
  functionals} of the $^*$-algebra $(C^\infty(M)[[\lambda]], \star,
\bar{\;})$ provide the physically relevant notion.  Recall that a
$\mathbb{C}[[\lambda]]$-linear functional $\omega:
C^\infty(M)[[\lambda]] \longrightarrow \mathbb{C}[[\lambda]]$ is
called positive if
\begin{equation}
    \label{eq:PosFun}
    \omega(\cc{f} \star f) \ge 0
\end{equation}
in the sense of formal power series, i.e. $a \in
\mathbb{R}[[\lambda]]$ with $a = \sum_{r=r_0}^\infty \lambda^r a_r$ is
positive if $a_{r_0} > 0$. Thus we enter the framework of
$^*$-algebras over $\ring{C} = \ring{R}(\im)$ with an ordered ring
$\ring{R}$, in our case $\ring{R} = \mathbb{R}[[\lambda]]$, as
indicated in the introduction.

Having a description for the states is still not enough for a quantum
mechanical theory: one also needs a way to describe
\emph{superpositions} of states as one of the most important physical
features of quantum physics. The naive convex combination of positive
linear functionals yields again a positive linear functional. However,
this does not correspond to superposition but to a mixed state.
Instead we need a linear structure for the states whence we have to
realize states as \emph{vector states} in a (pre-Hilbert) space where
the observable algebra acts on by a $^*$-representation. In general
this can not be accomplished for all states simultanously whence we
have to expect \emph{super-selection rules}. In particular, we have to
take into account \emph{all} representations of the observable algebra
in the beginning and find physical criteria for a selection of the
`interesting' ones afterwards. This way one is lead to the discussion
of the $^*$-representation theory of star product algebras. Of course,
this well-known line of argument applies to any sort of quantum theory
whose description is based on the observable algebra, like again the
algebraic quantum field theory.

%
%

\section{Notions of positivity and $^*$-representation theory}
\label{sec:positiv}

As the order structure of $\mathbb{R}[[\lambda]]$ plays a crucial role
in developing a representation theory for star products we shall now
recall some notions of positivity and $^*$-representation theory for
the framework of $^*$-algebras over $\ring{C} = \ring{R}(\im)$ where
$\ring{R}$ is an arbitrary ordered ring and $\im^2 = -1$. The
guideline for this development has been the rich theory of
$C^*$-algebras and many results and definitions can be seen as
algebraic analogs and generalizations of well-known results from
$C^*$-algebra theory, see
e.g.~\cite{landsman:1998a,raeburn.williams:1998a,lance:1995a}.

In fact, as we shall deal only with algebraic features, the more
appropriate analog is the theory of $^*$-algebras over $\mathbb{C}$
and their representations by (typically) \emph{unbounded} operators on
complex pre-Hilbert spaces: Representations of star product algebras
usually give differential operators, even in the cases with
convergence of the formal series, see e.g.
\cite{bordemann.neumaier.pflaum.waldmann:2003a}. Thus we shall take
the opportunity to relate our results to this framework, following
closely Schmüdgen's monography~\cite{schmuedgen:1990a}.

We start with a few basic concepts: A \emph{pre-Hilbert space} over
$\ring{C}$ is a $\ring{C}$-module $\mathcal{H}$ with sesquilinear
inner product $\SP{\cdot,\cdot}: \mathcal{H} \times \mathcal{H}
\longrightarrow \ring{C}$ such that $\SP{\phi, \psi} = \cc{\SP{\psi,
    \phi}}$ and $\SP{\phi, \phi} > 0$ for $\phi \ne 0$. A map $A:
\mathcal{H} \longrightarrow \mathcal{H}'$ between pre-Hilbert spaces
is called \emph{adjointable} if there exists a map $A^*: \mathcal{H}'
\longrightarrow \mathcal{H}$ with $\SP{\phi, A \psi} = \SP{A^*\phi,
  \psi}$ for all $\psi \in \mathcal{H}$ and $\phi \in \mathcal{H}'$.
In particular, $A$ is $\ring{C}$-linear and $A^*$ is unique, if it
exists at all. Then $\mathfrak{B}(\mathcal{H}, \mathcal{H}')$ denotes
the set of adjointable maps $\mathcal{H} \longrightarrow
\mathcal{H}'$. It follows that $\mathfrak{B}(\mathcal{H},
\mathcal{H}')$ is a $\ring{C}$-module, the map $A \mapsto A^*$ induces
a $\ring{C}$-antilineare involutive bijection
$\mathfrak{B}(\mathcal{H}, \mathcal{H}') \longrightarrow
\mathfrak{B}(\mathcal{H}',\mathcal{H})$ and the composition of
adjointable maps is again adjointable with $(AB)^* = B^*A^*$. This way
one obtains the category of pre-Hilbert spaces over $\ring{C}$ with
adjointable maps as morphisms. In particular,
$\mathfrak{B}(\mathcal{H}) = \mathfrak{B}(\mathcal{H}, \mathcal{H})$
is a $^*$-algebra over $\ring{C}$. Having $\mathfrak{B}(\mathcal{H})$
as `reference' $^*$-algebra one defines a \emph{$^*$-representation}
of an arbitrary $^*$-algebra $\mathcal{A}$ on a pre-Hilbert space
$\mathcal{H}$ to be a $^*$-homomorphism $\pi: \mathcal{A}
\longrightarrow \mathfrak{B}(\mathcal{H})$. Finally, an
\emph{intertwiner} $T$ between two $^*$-representations $(\mathcal{H},
\pi)$ and $(\mathcal{H}', \pi')$ is an adjointable map $T \in
\mathfrak{B}(\mathcal{H}, \mathcal{H}')$ with $T \pi(a) = \pi'(a) T$
for $a \in \mathcal{A}$. Two $^*$-representations are called
\emph{unitarily equivalent} if there exists a unitary intertwiner.
Comparing this with the case of $\ring{C} = \mathbb{C}$ we see that
the definition of $\mathfrak{B}(\mathcal{H})$ for a complex
pre-Hilbert space $\mathfrak{H}$ coincides with
$\mathcal{L}^+(\mathcal{H})$ in \cite[Prop.~2.1.8]{schmuedgen:1990a}
and the definition of a $^*$-representation is precisely the one in
\cite[Def.~2.1.13 and Def.~8.1.9]{schmuedgen:1990a}. On the other
hand, the intertwiners in \cite[Def.~8.2.1]{schmuedgen:1990a} are
required to be bounded, a notion which in our algebraic approach turns
out to be (unfortunately) rather useless, see \cite[App.~A and
B]{bordemann.waldmann:1998a} for examples. However, the difference
disappears for unitary equivalences as unitary maps on pre-Hilbert
spaces are of course bounded. In both cases one uses the intertwiners
as morphisms to define the category $\rep(\mathcal{A})$ of
$^*$-representations of $\mathcal{A}$.

Having Morita theory in mind, the above $^*$-representations are not
yet enough to obtain a clear formulation: we rather need more general
representation spaces where the inner products take their values in
general $^*$-algebras instead of $\ring{C}$ only. In order to
formulate the necessary positivity condition we have to specify
`positive elements' in a $^*$-algebra. In principle, there are many
possibilities: one needs to specify a \emph{m-admissible wedge}
$\mathcal{K} \subseteq \mathcal{A}$, i.e. a subset of Hermitian
elements closed under convex combinations, containing the elements of
the form $a^*a$ for $a \in \mathcal{A}$ and with $a^*\mathcal{K}a
\subseteq \mathcal{K}$ for all $a \in \mathcal{A}$, see
\cite[p.~22]{schmuedgen:1990a} for a definition in the case of
$\ring{C} = \mathbb{C}$, which immediately generalizes to arbitrary
$\ring{C}$.

For a $^*$-algebra there are \emph{two canonical} m-admissible
wedges. The first is given by
\begin{equation}
    \label{eq:Aplusplus}
    \mathcal{A}^{++} = \left\{a \in \mathcal{A} 
    \; \Big| \; a = \sum\nolimits_{i=1}^n \alpha_i a_i^*a_i
    \;
    \textrm{where}
    \;
    a_i \in \mathcal{A}
    \;
    \textrm{and}
    \;
    \ring{R} \ni \alpha_i > 0
    \right\},
\end{equation}
also denoted by $\mathcal{P}(\mathcal{A})$ in \cite{schmuedgen:1990a},
which is clearly the smallest m-admissible wedge and necessarily
contained in any other. The second is described as follows: recall
that a $\ring{C}$-linear functional $\omega: \mathcal{A}
\longrightarrow \ring{C}$ is called positive if $\omega(a^*a) \ge 0$
for all $a \in \mathcal{A}$. Then
\begin{equation}
    \label{eq:Aplus}
    \mathcal{A}^+ = \left\{a \in \mathcal{A}
    \; \big| \; \omega(a) \ge 0
    \;
    \textrm{for all positive}
    \; \omega
    \right\}
\end{equation}
is also a m-admissible wedge as $\omega_b(a) := \omega(b^*ab)$ is
again a positive functional for any $b \in \mathcal{A}$. In general
$\mathcal{A}^{++} \subsetneq \mathcal{A}^+$. Beside these two
canonical wedges one can obtain other m-admissible wedges by
additional choices or structures depending on which application one
has in mind.  One fairly general way to do this is to specify a subset
of positive functionals $S$ which is closed under convex combinations
and such that $\omega_b \in S$ for all $\omega \in S$ and $b \in
\mathcal{A}$. Then
\begin{equation}
    \label{eq:AplusS}
    \mathcal{A}_S^+ = \left\{a \in \mathcal{A}
    \; \big| \; \omega(a) \ge 0
    \;
    \textrm{for all}
    \; \omega \in S
    \right\}
\end{equation}
is again m-admissible and we have $\mathcal{A}_S^+ \subseteq
\mathcal{A}_{S'}^+$ for $S' \subseteq S$. In particular $\mathcal{A}^+
\subseteq \mathcal{A}^+_S$ for all such $S$. If $(\mathcal{H}, \pi)$
is a $^*$-representation then $S_\pi$ consisting of those positive
linear functionals, which are convex combinations of vector states in
$\pi$, i.e. $\omega_\varphi(a) = \SP{\varphi, \pi(a) \varphi}$ with
$\varphi \in \mathcal{H}$, turns out to be such a set $S$. In case of
an $O^*$-algebra $\mathcal{A}$, i.e. a $^*$-subalgebra of some
$\mathfrak{B}(\mathcal{H})$ this defines the \emph{positive cone} of
$\mathcal{A}$, using the defining representation of $\mathcal{A}$ on
$\mathcal{H}$, see \cite[Sect.~2.6]{schmuedgen:1990a}, leading in
particular to the notion of \emph{strong} positivity for
$O^*$-algebras. We note however, that only the two wedges
$\mathcal{A}^+$ and $\mathcal{A}^{++}$ are canonical and the others
depend on additional choices. Thus we shall stick to the choice
$\mathcal{A}^+$ for the `positive elements' of $\mathcal{A}$ and base
all our notions of positivity on $\mathcal{A}^+$. In particular, it
seems the `correct' choice in deformation quantization: At least for
the classical limit $C^\infty(M)$ this definition reproduces, by Riesz
representation theorem, the smooth functions on a manifold which have
non-negative values, see \cite[App.~B]{bursztyn.waldmann:2001a}.

%
%

\section{Pre-Hilbert modules and tensor products}
\label{sec:TensorProd}

In order to generalize the notion of $^*$-representations to more
general representation spaces we state the following definitions
\cite{bursztyn.waldmann:2003a:pre}. Let $\mathcal{D}$ be an auxilliary
$^*$-algebra over $\ring{C}$ and let $\HD$ be a right
$\mathcal{D}$-module. We always assume that modules have an underlying
compatible $\ring{C}$-module structure and that in the case of unital
$^*$-algebras the unit acts as identity on the module. Then a
\emph{$\mathcal{D}$-valued inner product} on $\mathcal{H}$ is a
$\ring{C}$-sesquilinear map (linear in the second argument)
\begin{equation}
    \label{eq:SPHDDD}
    \SPD{\cdot, \cdot} : \HD \times \HD \longrightarrow \mathcal{D}
\end{equation}
such that $\SPD{x,y} = (\SPD{y,x})^*$ and $\SP{x,y\cdot d} = \SPD{x,y}
d$ for all $x,y \in \HD$ and $d \in \mathcal{D}$. It is called
\emph{non-degenerate} if $\SPD{x,y} = 0$ for all $y$ implies $x=0$. In
this case $(\HD, \SPD{\cdot,\cdot})$ is called an inner product
$\mathcal{D}$-module. Then the characterization of adjointable
operators as in the case of pre-Hilbert spaces carries over and gives
$\ring{C}$-modules $\mathfrak{B}(\HD, \HDp)$ and a $^*$-algebra
$\mathfrak{B}(\HD)$. Note that adjointable maps are always
$\mathcal{D}$-linear to the right. The inner product is called
\emph{completely positive} if the matrix $(\SPD{x_i, x_j}) \in
M_n(\mathcal{D})^+$ is positive for all $n \in \mathbb{N}$ and $x_1,
\ldots, x_n \in \mathcal{D}$. If $\SPD{\cdot,\cdot}$ is non-degenerate
and completely positive then $(\HD, \SPD{\cdot,\cdot})$ is called a
\emph{pre-Hilbert $\mathcal{D}$-module}. Then a $^*$-representation
$\pi$ of a $^*$-algebra $\mathcal{A}$ on a pre-Hilbert
$\mathcal{D}$-module $\HD$ is a $^*$-homomorphism $\pi: \mathcal{A}
\longrightarrow \mathfrak{B}(\HD)$. In particular, $\HD$ becomes a
$(\mathcal{A}, \mathcal{D})$-bimodule. Clearly, the notion of
intertwiners carries over to this framework whence we obtain the
category $\rep[\mathcal{D}](\mathcal{A})$. If we drop the positivity
condition we still have the category $\smod[\mathcal{D}](\mathcal{A})$
of $^*$-representations of $\mathcal{A}$ on inner product
$\mathcal{D}$-modules.

We shall now show that a completely positive inner product
$\SPD{\cdot,\cdot}$ on $\HD$ can indeed be identified with a certain
completely positive map. To this end recall that the complex conjugate
$\mathcal{D}$-module $\DccH$ is defined to be $\cc{\mathcal{H}} =
\mathcal{H}$ as $\ring{R}$-module but now the $\ring{C}$- and
$\mathcal{D}$-module structure are defined by
\begin{equation}
    \label{eq:ccH}
    \alpha \cc{x} = \cc{x \cc{\alpha}}
    \quad
    \textrm{and}
    \quad
    d \cdot \cc{x} = \cc{x \cdot d^*},
\end{equation}
where $x \mapsto \cc{x}$ denotes the identity map $\mathcal{H}
\longrightarrow \cc{\mathcal{H}}$. This way, $\DccH$ becomes a left
$\mathcal{D}$-module. Then we consider the following
$(\mathcal{D},\mathcal{D})$-bimodule $\DED = \DccH \otimes_{\ring{C}}
\HD$. The map $I: \DED \longrightarrow \DED$ defined by
\begin{equation}
    \label{eq:IDED}
    I: \cc{x} \otimes y \mapsto \cc{y} \otimes x
\end{equation}
is a $\ring{C}$-antilinear $(\mathcal{D}, \mathcal{D})$-antibimodule
involution, i.e. we have $I(d \cdot z \cdot d') = {d'}^* \cdot I(z)
\cdot d^*$ for all $d,d' \in \mathcal{D}$ and $z \in \DED$. Moreover,
$\mathcal{H}^n$ becomes canonically a right $M_n(\mathcal{D})$-module
whence $M_n(\mathcal{E}) = \cc{\mathcal{H}^n} \otimes_{\ring{C}}
\mathcal{H}^n$ is a $(M_n(\mathcal{D}), M_n(\mathcal{D}))$-bimodule in
a canonical way. Clearly $I$ extends to $M_n(\mathcal{E})$ by
$I(z)_{ij} = I(z_{ji})$ for a matrix $z = (z_{ij}) \in
M_n(\mathcal{E})$. Still, we have $I(D \cdot z \cdot D') = (D')^*
\cdot I(z) \cdot D^*$ for $D, D' \in M_n(\mathcal{D})$ and $z \in
M_n(\mathcal{E})$. Now consider the subset
$\mathcal{K}_n(\mathcal{E})$ consisting of convex combinations of
elements $z \in M_n(\mathcal{E})$ with $z_{ij} = \cc{x_i} \otimes x_j$
where $x_1, \ldots, x_n \in \mathcal{H}$. Clearly
$I(\mathcal{K}_n(\mathcal{E})) = \mathcal{K}_n(\mathcal{E})$ and $D^*
\mathcal{K}_n(\mathcal{E}) D \subseteq \mathcal{K}_n(\mathcal{E})$ for
all $D \in M_n(\mathcal{D})$. This way, $\mathcal{E}$ becomes a
\emph{matrix-ordered space} in the sense of
\cite[Def.~11.1.1]{schmuedgen:1990a}, even in the sense of
$(\mathcal{D}, \mathcal{D})$-bimodules, i.e. $\ring{C}$ replaced by
$\mathcal{D}$ everywhere. This allows to speak of a \emph{completely
  positive map} $\Phi: \mathcal{E} \longrightarrow \mathcal{D}$, i.e.
a $(\mathcal{D}, \mathcal{D})$-bimodule map such that $\Phi(z)^* =
\Phi(I(z))$ and
\begin{equation}
    \label{eq:PhiComPos}
    \Phi\left(\mathcal{K}_n(\mathcal{E})\right)
    \subseteq M_n(\mathcal{D})^+
\end{equation}
for all $n \in \mathbb{N}$. Again, we have to extend the definition of
completely positive maps as in \cite[Def.~11.1.3]{schmuedgen:1990a} to
the case of $(\mathcal{D}, \mathcal{D})$-bimodules.
\begin{proposition}
    \label{proposition:ComPos}
    Let $\HD$ be a right $\mathcal{D}$-module. Then a
    $\mathcal{D}$-valued inner product is equivalent to a
    $(\mathcal{D}, \mathcal{D})$-bimodule map $\Phi: \DED = \DccH
    \otimes_{\ring{C}} \HD \longrightarrow \mathcal{D}$ with
    $\Phi(z)^* = \Phi(I(z))$ via
    \begin{equation}
        \label{eq:PhiisSPD}
        \SPD{x,y} = \Phi(\cc{x} \otimes y).
    \end{equation}
    Moreover, $\SPD{\cdot,\cdot}$ is completely positive iff $\Phi$ is
    completely positive.
\end{proposition}
\begin{proof}
    This is now a straightforward verification.
\end{proof}

We turn now to tensor products. Consider $\BEA \in
\smod[\mathcal{A}](\mathcal{B})$ and $\CFB \in
\smod[\mathcal{B}](\mathcal{C})$. Then on $\CFB \tensor[B]
\BEA$, which is a $(\mathcal{C}, \mathcal{A})$-bimodule, one has a
$\mathcal{A}$-valued inner product defined by
\begin{equation}
    \label{eq:RieffelSP}
    \SPFEA{x_1 \otimes y_1, x_2 \otimes y_2}
    =
    \SPEA{y_1, \SPFB{x_1,x_2} \cdot y_2}.
\end{equation}
As $\SPFEA{\cdot,\cdot}$ may be degenerate one needs to divide by the
degeneracy space which is possible as it is invariant under the
$(\mathcal{C},\mathcal{A})$-bimodule structure.  As quotient one
obtains an inner product $\mathcal{A}$-module $\CFB \tensM[B] \BEA$
with a $^*$-representation of $\mathcal{C}$. Hence one ends up with a
functor
\begin{equation}
    \label{eq:TensHat}
    \tensM[B]: 
    \smod[\mathcal{B}](\mathcal{C}) 
    \times
    \smod[\mathcal{A}](\mathcal{B})
    \longrightarrow
    \smod[\mathcal{A}](\mathcal{C}).
\end{equation}
Clearly, $\tensM$ is associative up to the usual canonical
isomorphisms. It turns out that completely positive inner products
behave well under this tensor product
\cite[Thm.~4.7]{bursztyn.waldmann:2003a:pre} whence $\tensM[B]$
restricts to a functor
\begin{equation}
    \label{eq:tensMPos}
    \tensM[B]: 
    \rep[\mathcal{B}](\mathcal{C}) 
    \times
    \rep[\mathcal{A}](\mathcal{B})
    \longrightarrow
    \rep[\mathcal{A}](\mathcal{C}).
\end{equation}

By fixing e.g. the first argument of $\tensM[B]$ one arrives at the
algebraic analog of Rieffel's induction procedure for
$^*$-representations of $C^*$-algebras
\cite{rieffel:1974a,rieffel:1974b}, see also
\cite[Ex.~4.9]{bursztyn.waldmann:2003a:pre}, which itself can be seen
as a generalization \cite[Prop.~4.7]{bursztyn.waldmann:2001a} of the
well-known \emph{GNS construction} of $^*$-representations out of a
positive linear functional, see e.g. any textbook on $C^*$-algebras
and \cite[Sect~8.6]{schmuedgen:1990a} for a detailed exposition in the
context of $^*$-algebras over $\mathbb{C}$ and
\cite{bordemann.waldmann:1998a} for applications to deformation
quantization.

%
%

\section{Strong Morita equivalence and the Picard groupoid}
\label{sec:Morita}

Morita equivalence was first developed by Morita~\cite{morita:1958a}
in the purely algebraic framework of associative rings, see
e.g.~\cite{bass:1968a}. Later, Rieffel transfered these ideas to
$C^*$/$W^*$-algebras leading to the notion of strong Morita
equivalence \cite{rieffel:1974a,rieffel:1974b}. Ara specialised the
general notion to rings with involution coining the notion of
$^*$-Morita equivalence \cite{ara:1999b,ara:1999a} and
Bursztyn-Waldmann discussed the case of $^*$-algebras over $\ring{C} =
\ring{R}(\im)$, generalizing Rieffel's strong Morita equivalence to
this situation, therefor also called strong Morita equivalence
\cite{bursztyn.waldmann:2003a:pre,bursztyn.waldmann:2001a}. A
reasonable notion of Morita equivalence in the last class of algebras
in the list in the introduction, the $^*$-algebras over $\mathbb{C}$,
seems still missing. Here one wants to go beyond the purely algebraic
treatment taking into account the much more refined notions of
$^*$-representation theory available here \cite{schmuedgen:1990a}. In
particular, given some locally convex topologies on the $^*$-algebras,
one would like to have some notion of Morita equivalence respecting
these extra structures.

For $C^*$-algebras all the possible notions of Morita equivalence give
the same equivalence relation, see
\cite{beer:1982a,ara:1999b,bursztyn.waldmann:2001b} while usually one
only has the implications
\begin{equation}
    \label{eq:MEImplication}
    \textrm{Strong Morita equivalence}
    \Rightarrow
    \textrm{$^*$-Morita equivalence}
    \Rightarrow
    \textrm{Morita equivalence}.
\end{equation}
Star product algebras are shown to behave like $C^*$-algebras
concerning \eqref{eq:MEImplication}: here also the ring-theoretical
Morita equivalence implies strong Morita equivalence
\cite{bursztyn.waldmann:2002a}. However, the situation becomes more
complicated if one asks in how many (essentially different) ways two
$^*$-algebras can be ($^*$- resp. strongly) Morita equivalent. This
information is encoded in the ($^*$- resp. strong) Picard groupoid
which we shall now describe.

Thus we first have to recall the definition of strong and $^*$-Morita
equivalence in some more detail. Let $\BEA$ be a $(\mathcal{B},
\mathcal{A})$-bimodule endowed with a $\mathcal{A}$-valued inner
product $\SPA{\cdot,\cdot}$ and a $\mathcal{B}$-valued inner product
$\BSP{\cdot,\cdot}$. As $\mathcal{B}$ acts on $\mathcal{E}$ from the
left the $\mathcal{B}$-valued inner product is linear (to the left) in
the first argument. An inner product $\SPA{\cdot,\cdot}$ is called
\emph{full} if the span of all $\SPA{x,y}$ in $\mathcal{A}$ coincides
with $\mathcal{A}$. This is equivalent to say that the map $\Phi$
associated to $\SPA{\cdot,\cdot}$ as in
Proposition~\ref{proposition:ComPos} is surjective. The $(\mathcal{B},
\mathcal{A})$-bimodule $\BEA$ with inner products $\BSP{\cdot,\cdot}$
and $\SPA{\cdot,\cdot}$ is called a \emph{$^*$-equivalence bimodule} if
$\mathcal{B} \cdot \mathcal{E} = \mathcal{E} = \mathcal{E} \cdot
\mathcal{A}$, both inner products are non-degenerate and full, and
they satisfy
\begin{equation}
    \label{eq:Compatible}
    \SPA{x, b \cdot y} = \SPA{b^* \cdot x, y},
    \quad
    \BSP{x \cdot a, y} = \BSP{x, y \cdot a^*},
    \quad
    \textrm{and}
    \quad
    \BSP{x,y} \cdot z = x \cdot \SPA{y,z}
\end{equation}
for all $x,y,z \in \mathcal{E}$, $b \in \mathcal{B}$ and $a \in
\mathcal{A}$. If in addition the inner products are completely
positive then $\BEA$ is called a \emph{strong equivalence bimodule},
see~\cite{ara:1999a,bursztyn.waldmann:2003a:pre}. If such a bimodule
exists the algebras are called \emph{$^*$-Morita equivalent} or
\emph{strongly Morita equivalent}, respectively. It turns out that the
$\tensM$-tensor product of equivalence bimodules is again an
equivalence bimodule. This leads to the statements that $^*$-Morita
equivalence as well as strong Morita equivalence is indeed an
equivalence relation within the class of $^*$-algebras which are
non-degenerate and idempotent. The condition of non-degeneracy and
idempotency has to be put since otherwise both notions fail to be
reflexive. In particular, unital $^*$-algebras are always
non-degenerate and idempotent.

The \emph{strong Picard groupoid} is now constructed as follows.
$\StrPic(\mathcal{B}, \mathcal{A})$ consists of (isometric)
isomorphism classes of equivalence bimodules $\BEA$ and
$\StrPic(\mathcal{B}, \mathcal{A})$ is viewed as space of arrows of a
(large) groupoid with class of units given by the $^*$-algebras
themselves. The composition law
\begin{equation}
    \label{eq:CompositionStrPic}
    \StrPic(\mathcal{C}, \mathcal{B}) \times
    \StrPic(\mathcal{B}, \mathcal{A})
    \longrightarrow
    \StrPic(\mathcal{C}, \mathcal{A})
\end{equation}
is given by the tensor product $\tensM$ and the units are the trivial
bimodules $\AAA$ with the canonical inner products $\SPA{a,b} = a^*b$
and $\ASP{a,b} = ab^*$. Analogously one defines the $^*$-Picard
groupoid $\starPic(\cdot,\cdot)$ based on $^*$-equivalence bimodules.
It turns out that this defines indeed a groupoid
$\StrPic(\cdot,\cdot)$, see
\cite[Thm.~6.1]{bursztyn.waldmann:2003a:pre} for an alternative
construction and e.g.~\cite{bass:1968a,benabou:1967a} for the
ring-theoretical Picard groupoid. Finally, one observes that Rieffel's
original version of strong Morita equivalence for $C^*$-algebras
(which involves additional completeness properties) gives a
$C^*$-algebraic strong Picard groupoid $\StrPic_{C^*}$, see
\cite{brown.green.rieffel:1977a}, which turns out to be equivalent to
the considerably easier strong Picard groupoid $\StrPic$, restricted
to the corresponding Pedersen ideals, see
\cite[Cor.~6.11]{bursztyn.waldmann:2003a:pre}. Thus $\StrPic$ encodes
the complete strong Morita theory for $C^*$-algebras.

In general, the understanding of Mortia theory is equivalent to the
understanding of the Picard groupoid in each of its flavours
above. There are two basic questions:
\begin{enumerate}
\item Which $^*$-algebras $\mathcal{A}$, $\mathcal{B}$ are strongly
    Morita equivalent, i.e. for which is $\StrPic(\mathcal{B},
    \mathcal{A}) \ne \emptyset$?
\item How many different self-equivalences does $\mathcal{A}$ have,
    i.e. what is the \emph{strong Picard group} $\StrPic(\mathcal{A}) =
    \StrPic(\mathcal{A},\mathcal{A})$ of $\mathcal{A}$?
\end{enumerate}
Konwing this we have a complete description of $\StrPic(\cdot,\cdot)$
as it is a groupoid and a groupoid is determined by its orbits and
isotropy groups. In particular, if $\StrPic(\mathcal{B}, \mathcal{A})
\ne \emptyset$ then $\StrPic(\mathcal{B}, \mathcal{A}) \cong
\StrPic(\mathcal{A})$ as sets since isotropy groups (here
$\StrPic(\mathcal{A})$) act transitively and freely on the spaces of
arrows. Hence in this situation $\StrPic(\mathcal{A}) \cong
\StrPic(\mathcal{B})$ as groups where each element of
$\StrPic(\mathcal{B},\mathcal{A})$ provides an isomorphism.

Another type of question is how the different flavours of Picard
groupoids for $^*$-algebras over $\ring{C}$ are related. It is rather
easy to see that `forgetting' the positivity of the inner products or
`forgetting' the inner products at all yields groupoid morphisms
\begin{equation}
    \label{eq:PicPicPic}
    \StrPic \longrightarrow \starPic,
    \quad
    \StrPic \longrightarrow \Pic,
    \quad
    \textrm{and}
    \quad
    \starPic \longrightarrow \Pic
\end{equation}
which are in general neither injective nor surjective. The morphism
$\StrPic \longrightarrow \starPic$ usually fails to be surjective as
there may be also other inner products on equivalence bimodules having
different `signatures'. The morphism $\starPic \longrightarrow \Pic$
fails to be injective for the same reason. However, also $\StrPic
\longrightarrow \Pic$ shows a rather rich and complicated behaviour,
even for very `nice' classes of $^*$-algebras like $C^*$-algebras. For
$C^*$-algebras as well as for star product algebras one can prove
injectivity but surjectivity fails in general. The defect of being
surjective can be described rather explicit in terms of certain
non-inner automorphisms, see
\cite[Sect.~7,8]{bursztyn.waldmann:2003a:pre} for a detailed analysis
and examples.

The last type of question we want to mention is about formal
deformations of $^*$-algebras. Thus assume that
$\boldsymbol{\mathcal{A}} = (\mathcal{A}[[\lambda]], \star)$ and
$\boldsymbol{\mathcal{B}} = (\mathcal{B}[[\lambda]], \star')$ are
Hermitian formal deformations of $\mathcal{A}$, $\mathcal{B}$. Then
one wants to relate $\StrPic(\boldsymbol{\mathcal{B}},
\boldsymbol{\mathcal{A}})$ and $\StrPic(\mathcal{B}, \mathcal{A})$.
Here one heavily uses the fact that formal deformations and ordered
rings fits together nicely as $\ring{R}[[\lambda]]$ is canonically
ordered again. While for $\Pic(\cdot, \cdot)$ and $\starPic(\cdot,
\cdot)$ it is rather easy to see that the `classical limit' map $\cl$,
which sets $\lambda = 0$, induces a groupoid morphism $\cl_*$. This
turns out to be more complicated in the case of
$\StrPic(\cdot,\cdot)$. Here it is only true for \emph{completely
  positive deformations}. Recall that a Hermitian deformation
$\boldsymbol{\mathcal{A}}$ is called completely positive if for every
$\ring{C}$-linear positive functional $\Omega_0: M_n(\mathcal{A})
\longrightarrow \ring{C}$ one finds a deformation into a
$\ring{C}[[\lambda]]$-linear positive functional $\Omega = \Omega_0 +
o(\lambda): M_n(\boldsymbol{\mathcal{A}}) \longrightarrow
\ring{C}[[\lambda]]$. Luckily, Hermitian star products are always
completely positive deformations
\cite[Thm.~8.5]{bursztyn.waldmann:2003a:pre} whence we have a groupoid
morphism
\begin{equation}
    \label{eq:clLimit}
    \cl_*: \StrPic(\boldsymbol{\mathcal{B}}, \boldsymbol{\mathcal{A}})
    \longrightarrow
    \StrPic(\mathcal{B}, \mathcal{A})
\end{equation}
for star product algebras and, more generally, completely positive
deformations. This actually follows from the considerations in
\cite[Lem.~9.3]{bursztyn.waldmann:2001a}. Again, for all three
versions $\Pic$, $\starPic$ and $\StrPic$ the classical limit $\cl_*$
shows a rich and non-trivial behaviour, see
\cite{bursztyn.waldmann:2002a:pre} for a detailed discussion of the
case of $\Pic$ and
\cite{jurco.schupp.wess:2002a,bursztyn.waldmann:2002a} for a
classification of (strongly) Morita equivalent star products resulting
from these considerations concerning the classical limit map.

%
%

\section{The Picard groupoid in action}
\label{sec:Action}

In this last section we illustrate how one can obtain Morita
invariants of algebras by looking at `actions' of the Picard
groupoids. We shall not give an axiomatic theory for actions of a large
groupoid (or even bigroupoid) but mention just three examples, the
results of which are classical:
\begin{enumerate}
\item As already mentioned $\StrPic(\mathcal{A})$ is an invariant
    under strong Morita equivalence. This can be seen as result of the
    groupoid action of $\StrPic$
    \begin{equation}
        \label{eq:PiconItself}
        \StrPic(\mathcal{B}, \mathcal{A}) \times \StrPic(\mathcal{A})
        \longrightarrow
        \StrPic(\mathcal{B})
    \end{equation}
    on its isotropy groups by simply using the left multiplications.
    
\item Let $\mathcal{A}$, $\mathcal{B}$ be unital. Then the Hermitian
    $K_0$-group of $\mathcal{A}$ is the Grothendieck group
    $K^H_0(\mathcal{A})$ of the semi-group of isometric isomorphism
    classes of finitely generated projective right
    $\mathcal{A}$-modules with strongly non-degenerate completely
    positive inner products. Here strongly non-degenerate means that
    $x \in \mathcal{E} \mapsto \SPA{x, \cdot} \in \mathcal{E}^* =
    \Hom_{\mathcal{A}}(\mathcal{E}, \mathcal{A})$ is a bijection. Now
    let $\BEA$ be a strong equivalence bimodule and let $\HB$ be a
    finitely generated projective pre-Hilbert module with strongly
    non-degenerate inner product. Then $\HpA = \HB \tensM[B] \BEA$ is
    again finitely generated and projective as right
    $\mathcal{A}$-module and the completely positive inner product is
    again strongly non-degenerate. Passing to isometric isomorphism
    classes this yields a `groupoid action'
    \begin{equation}
        \label{eq:KBPicBAKA}
        K^H_0(\mathcal{B}) \times 
        \StrPic(\mathcal{B}, \mathcal{A})
        \longrightarrow
        K^H_0(\mathcal{A})
    \end{equation}
    of the strong Picard groupoid on the Hermitian $K_0$-groups of
    unital $^*$-algebras over $\ring{C}$. From this we can see two
    things: first $K_0(\mathcal{A})$ carries a natural
    $\StrPic(\mathcal{A})$-representation, as the tensor product
    $\tensM[\mathcal{A}]$ is clearly compatible with the direct
    orthogonal sum, and second $K_0^H(\mathcal{A}) \cong
    K_0^H(\mathcal{B})$ as $\StrPic(\mathcal{A}) \cong
    \StrPic(\mathcal{B})$ representation spaces. This observation was
    brought to our attention by Alan Weinstein.

\item Let $\mathcal{D}$ be an auxilliary $^*$-algebra and denote by
    $\Rep[\mathcal{D}](\mathcal{A})$ those $^*$-representations
    $(\mathcal{H}, \pi) \in \rep[\mathcal{D}](\mathcal{A})$ where
    $\pi(\mathcal{A})\mathcal{H} = \mathcal{H}$. Then given a strong
    equivalence bimodule $\BEA$ and $(\mathcal{H}, \pi) \in
    \Rep[\mathcal{D}](\mathcal{A})$ we obtain $\BEA \tensM[A]
    \mathcal{H} \in \Rep[\mathcal{D}](\mathcal{B})$ and the fucntor
    $\mathcal{H} \mapsto \BEA \tensM[A] \mathcal{H}$ implements an
    equivalence of the categories $\Rep[\mathcal{D}](\mathcal{A})$ and
    $\Rep[\mathcal{D}](\mathcal{B})$, see
    \cite[Cor.~5.16]{bursztyn.waldmann:2003a:pre}. Note that this was
    one of the original motivations for the definition of strong
    Morita equivalence. Reinterpreting this Morita invariant we see
    that this comes from an `action'
    \begin{equation}
        \label{eq:PicRepRep}
        \StrPic(\mathcal{B}, \mathcal{A}) 
        \times
        \Rep[\mathcal{D}](\mathcal{A})
        \longrightarrow
        \Rep[\mathcal{D}](\mathcal{B}),
    \end{equation}
    which, however, is only defined up to natural isomorphisms and
    unitary equivalences. Thus the probably better framework would be
    to consider the Picard \emph{bigroupoid} instead, i.e not
    identifying bimodules up to isomorphism, see
    \cite{benabou:1967a}. Then \eqref{eq:PicRepRep} would better read
    as an `action' of the Picard bigroupoid on the categories of
    $^*$-representations.
\end{enumerate}

Let us finally mention that the above examples have their well-known
counterparts for $\StrPic$ replaced by $\Pic$ or $\starPic$ leading
the known Morita invariants in these contexts. Beside the above
examples there are other Morita invariants which can be obtained by
actions of the Picard groupoid like the moduli spaces of formal
(Hermitian) deformations $\Def(\mathcal{A})$, see the discussion in
\cite[Sect.~3.3]{bursztyn.waldmann:2002a:pre} and
\cite{bursztyn:2002a}.

%
%

\begin{footnotesize}
    \renewcommand{\arraystretch}{0.5} 

\end{footnotesize}

\end{document}